\documentclass[11pt]{article}
\usepackage{amsmath,amsthm,amssymb,graphicx}
\usepackage[margin=1in]{geometry}
\usepackage[colorlinks=true,linkcolor=black,citecolor=black,urlcolor=blue]{hyperref}
\newcommand{\doi}[1]{\href{https://doi.org/#1}{\nolinkurl{doi:#1}}}

\title{On the Role of Fisher Information, Variance and Asymptotics in Statistical Inference}
\author{Paul W. Vos\thanks{Department of Public Health, East Carolina University,
Greenville, North Carolina. Email: \texttt{vosp@ecu.edu}. ORCID: 0000-0001-9996-5627.}}
\date{November 2026}
\begin{document}
\maketitle

\begin{abstract}
Bahadur's \emph{Lectures on the Theory of Estimation} close with a chapter on asymptotics: asymptotic normality of the maximum likelihood estimator, Fisher's bound on its asymptotic variance with Hodges's counterexample, and local asymptotic normality. A companion paper read the earlier lectures with one modification, replacing the point estimator by a generalized estimator---Fisher's continuum of significance tests---assessed by information rather than variance. This paper reads Chapter 8 the same way. At each of its three parts the classical statement is an approximation about a point estimator, while the corresponding statement about the score is exact at every sample size. Superefficiency does not arise, because information is bounded by the Fisher information at every $n$; Hodges's estimator, computed exactly, is fully efficient at the origin and loses information only where its mean square error is worst. The information of a generalized estimator is argued to be Fisher information: the mean rate at which a continuum of tests responds to the parameter, standardized by its variance. The resulting Fisher upper bound differs geometrically, not merely algebraically, from the Cram\'er--Rao lower bound. The extension is not empty: the $t$-test continuum loses about half an observation's worth of the information in $\bar X$ at every $n$.
\end{abstract}

\section{Introduction}
\label{sec:intro}

Three roles are commonly assigned to asymptotic theory in statistical inference. It supplies a sampling distribution for an estimator whose exact distribution is unavailable; it supplies, through that distribution, a variance by which the estimator is assessed; and it supplies, through the Fisher information, a bound against which the variance is judged. The three roles are not on the same footing. The first is delivered by the central limit theorem. The second is not: convergence in distribution says nothing about convergence of moments, and the variance of the limit is in general not the limit of the variances. The third rests on the second and inherits its defect. Bahadur, in the closing chapter of his \emph{Lectures on the Theory of Estimation} \cite{Bahadur}, records the point with his usual care: immediately after stating that $\sqrt n(\hat\theta_n-\theta)\to N(0,1/I_1(\theta))$ he adds that the theorem ``does not assert that $E_\theta(\hat\theta_n)=\theta+o(1)$ or that $\operatorname{Var}_\theta(\hat\theta_n)=1/(nI_1(\theta))+o(1/n)$.'' The theorem is used, in practice, for exactly what it does not assert.

This paper argues that the difficulty is not with asymptotics but with the object asked to carry it. The classical theory assesses a point estimator---a function on the sample space, taking values in the parameter space---and describes it by moments, which are what the central limit theorem does not deliver. A generalized estimator \cite{VosWu,NewLook} is a function on the sample space and the parameter space jointly: for each sample it is a curve over $\Theta$, read as Fisher's continuum of significance tests, one for each hypothesis $\theta=\theta_0$. It is assessed by its information, the squared mean rate at which the standardized curve responds as $\theta$ moves. Information is defined on the object the central limit theorem is about, the score, and it is exact at every sample size; variance is defined on the object the theorem is used for, the estimator, and it is not. Once the object is changed, the three roles separate: the first remains, for the shape of a sampling distribution when it must be approximated; the second and third are not needed, because the assessment is exact without them.

The paper is the third of a series applying \cite{VosWu} to questions the statistical community has framed in classical terms. The first \cite{MSE} took up mean square error and the James--Stein paradox; the second \cite{NewLook} read Lectures 5--27 of Bahadur's notes, making one small modification at each place Bahadur marks a strain in the classical theory, and left Chapter 8, the asymptotic chapter, aside. This paper reads Chapter 8, whose content is the standard asymptotic theory of estimation as it appears in the textbooks \cite{CoxHinkley1974,LehmannCasella1998,CasellaBerger2024,vanderVaart1998}. It is short---eight pages, Lectures 28--30---and has three parts: the asymptotic normality of the maximum likelihood estimator; Fisher's bound on the asymptotic variance, Hodges's superefficient estimator, and the theorems of Le Cam and H\'ajek that recover the bound; and the local asymptotic normality condition. Section~\ref{sec:ch8} takes the parts in order. At each, the classical statement is an approximation about a point estimator and the corresponding statement about the score is exact; at each, the classical apparatus (regularity, measure-zero exceptions, the convolution theorem) is an asymptotic counterpart of the finite-sample apparatus catalogued in \cite{NewLook}---unbiasedness, admissibility, minimaxity---and answers the same need. Hodges's estimator is computed exactly rather than asymptotically, and under the information criterion there is nothing super about it.

Section~\ref{sec:fisher} motivates an alternative definition of the Fisher information in a statistic, one that makes the idea more intuitive and that extends at once from statistics to continua of tests. In \cite{VosWu} the information of a generalized estimator, written $\Lambda(\tau)$, was motivated by geometry; here the motivation is by statistical principles. Fisher information is defined by the two ingredients Fisher's own construction requires---the mean derivative of a continuum of tests in the parameter, which makes the assessment dependent on the family, and the variance, which standardizes it---and so defined it is $\Lambda(\tau)$. For the score the Bartlett identity makes the two ingredients numerically equal, and Fisher information is then usually defined by the variance alone; the collapse is a convenience for the score and an obstacle for the concept. The main result of \cite{VosWu} is restated with every ingredient visible as the Fisher upper bound, $\Lambda(\tau)\le I_M$, and set beside the Cram\'er--Rao lower bound, from which it differs by a rearrangement algebraically and fundamentally in its geometry. With the general definition in hand a question that has no classical formulation becomes immediate, and its answer for the $t$-test continuum is exact and, we think, new: Student's $t$-test loses about half an observation's worth of the information in $\bar X$ at every $n$, a loss first-order asymptotics cannot see.

What is new here: the reading of Chapter 8 (Section~\ref{sec:ch8}); the exact finite-sample computation for Hodges's estimator (Section~\ref{sec:hodges-finite}); the two-ingredient account of Fisher information, the proposal that $\Lambda$ be called Fisher information, and the geometric contrast with Cram\'er--Rao (Sections~\ref{sec:two}--\ref{sec:metric}); and the efficiency of the $t$-test continuum (Section~\ref{sec:t}). Recalled from \cite{VosWu,NewLook,MSE} and used without new argument: the definition of a generalized estimator and of its information, the bound with the score attaining it, the two maps from point estimators to generalized estimators, and the finite-sample counterpart of Efron's curvature. The parameter is scalar except where Bahadur's own statement is not: the local asymptotic normality condition of Lecture 30 is stated for $\theta\in\mathbb{R}^p$ with $a\in\mathbb{R}^p$ a direction, and we read it as he states it, drawing on the multiparameter theory of \cite[Sec.~3]{VosWu} only for the reading; nuisance parameters enter only in the $t$ example, through the orthogonalization of \cite{VosWu}. Section~\ref{sec:ch8} begins by recalling the notation of \cite{NewLook}; Section~\ref{sec:discussion} closes with what has and has not been shown.

\section{Bahadur's Chapter 8, read with generalized estimators}
\label{sec:ch8}

Chapter 8 of the lectures occupies eight pages (Lectures 28--30) and has three parts: the asymptotic normality of the maximum likelihood estimator, Fisher's bound on the asymptotic variance with Hodges's counterexample and the theorems of Le Cam and H\'ajek that repair it, and the local asymptotic normality condition under which the whole account extends beyond the iid case. We take the parts in order. At each, the classical statement is an approximation about a point estimator, and at each the corresponding statement about the score is exact at every sample size. The pattern is the one already seen in the odds ratio of \cite[Sec.~5.4]{NewLook}, where Bahadur's remedy for a plug-in estimate with infinite expectation was ``asymptotic---for large $n$ the offending term is small with high probability.'' Chapter 8 is that remedy applied systematically.

\paragraph{Notation recalled.}
We use Bahadur's notation as in \cite{NewLook}. A sample $s$ (in this chapter, $s_n=(x_1,\dots,x_n)$) takes values in a sample space $S$; the family is $\{P_\theta:\theta\in\Theta\}$ with $\Theta$ an open interval, densities $\ell_\theta=\ell(\theta\mid s)$ with respect to a fixed measure, and $L_\theta=\log\ell_\theta$. The likelihood ratio of $\delta$ to $\theta$ is $\Omega_{\delta,\theta}=\ell_\delta/\ell_\theta$, and the score at $\theta$ is
\[
\gamma^{(1)}_\theta(s)=\frac{\partial}{\partial\theta}L(\theta\mid s)=\frac{\partial}{\partial\delta}\Omega_{\delta,\theta}\Big|_{\delta=\theta},
\]
with $E_\theta\gamma^{(1)}_\theta=0$ and Fisher information $I=E_\theta[(\gamma^{(1)}_\theta)^2]$; for an iid sample $I=nI_1$ with $I_1$ the information in one observation. $V_\theta=L^2(S,P_\theta)$ is the Hilbert space of square-integrable functions under $P_\theta$, and Bahadur's estimate is a function $t$ on $S$ lying in every $V_\theta$.

A \emph{generalized estimator} \cite[Def.~2]{NewLook} is a function $\tau:S\times\Theta\to\mathbb{R}$ such that (i) $\tau(\cdot,\theta)\in V_\theta$ with $E_\theta\,\tau(s,\theta)=0$ for every $\theta$; (ii) $\tau(s,\cdot)$ is continuously differentiable on $\Theta$ for almost every $s$; (iii) $\operatorname{Var}(\tau)=E_\theta[\tau(s,\theta)^2]>0$ for every $\theta$. For fixed $s$ the estimate $\tau_s=\tau(s,\cdot)$ is a curve over $\Theta$, read as a continuum of test statistics: $\tau_s(\theta_0)$ is the observed value of a statistic that has mean zero under $P_{\theta_0}$. The score is the archetype. Its information is
\[
\Lambda(\tau)=\frac{\big(E[\partial\tau/\partial\theta]\big)^2}{\operatorname{Var}(\tau)},
\]
a function on $\Theta$, with $\Lambda(\tau)\le I$ for every $\tau$ and $\Lambda(\gamma^{(1)})=I$ \cite[Thm.~1]{NewLook}; $\mathrm{Eff}(\tau)=\Lambda(\tau)/I$. A bar denotes standardization, $\bar\tau=\tau/\sqrt{\operatorname{Var}(\tau)}$, so that the standardized score is $\bar\gamma^{(1)}=\gamma^{(1)}/\sqrt{I}$. When a particular estimator is being assessed we also write $I_\tau=\Lambda(\tau)$, in the notation Fisher used for the information $I_t$ in a statistic $t$; Section~\ref{sec:fisher} argues that this is more than an analogy. The functional notation $\Lambda(\cdot)$ is kept for statements, such as the bound, that concern information as a function on the space of generalized estimators. Expectations and variances of a generalized estimator carry no subscript, the second argument naming the distribution.

A point estimator $t$ induces two generalized estimators. The \emph{expectation map} gives $\tau_{E(t)}(s,\theta)=t(s)-E_\theta(t)$, built from the mean function of $t$. The \emph{score map} gives $\tau_t$, the score of the marginal family $\{P_\theta\circ t^{-1}\}$ of the statistic $t$, composed with $t$ so as to be a function on $S\times\Theta$; it depends on $t$ only through the $\sigma$-field $t$ generates, and $\Lambda(\tau_t)=I_t$, the Fisher information in $t$ \cite[eq.~(13)]{NewLook}.

In the other direction, a generalized estimator induces a point estimator through its roots. The graph of the estimate $\tau_s$ is a curve in $\Theta\times\mathbb{R}$, and the point estimate is where the curve crosses the horizontal axis: $\hat\theta(s)$ solves $\tau_s(\theta)=0$. For the score this is the maximum likelihood estimate, and for $\tau_{E(t)}$ with $t$ unbiased it is $t(s)$ itself. The root need not exist---the binomial score at $y=0$ is $-n/(1-p)$, negative on all of $(0,1)$---and need not be unique; the generalized estimate is defined and usable in either case, which is the existence payoff of \cite[Sec.~2.2]{NewLook}.

The distinction between point and generalized estimators can be seen in two distributions \cite[Example~1]{VosWu}, and the distinction recurs below. Fixing $\theta$ and letting $s$ vary gives the \emph{vertical distribution at $\theta$}: the values $\tau(s,\theta)$ of the curves on the vertical line over $\theta$, a distribution on $\mathbb{R}$ that has mean zero and, once standardized, variance one, in every parameterization. Letting the curves cross the horizontal line at height zero gives the \emph{horizontal distribution at $\theta$}: the roots $\{\theta':\tau_s(\theta')=0\}$ when $s$ is distributed according to $P_\theta$, a distribution on $\Theta$---the sampling distribution of the point estimate when the root exists---whose mean and variance depend on the parameterization and need not exist. As $\theta$ varies each kind gives a continuum of distributions; we say continuum rather than family because $\theta$ is more than an index---it labels the points of a family that already has a smooth structure, and the continuum inherits that structure. The variance of a point estimator, the classical assessment, is computed on horizontal distributions; Fisher information is computed on vertical ones.

\subsection{Lecture 28: asymptotic normality, and a note}
\label{sec:L28}

Bahadur's setting is the regular case: $\Theta$ an open interval, $\ell(\theta\mid x)>0$, $L(\theta\mid x)=\log \ell(\theta\mid x)$ with two continuous derivatives, $E_\theta L'(\theta\mid x)=0$, and a domination condition on $L''$ in a neighborhood of each $\theta$ (his Condition $(\ast)$): for each $\theta$ there is an $\varepsilon>0$ such that $\max_{|\delta-\theta|\le\varepsilon}|L''(\delta\mid x)|$ has finite expectation under $P_\theta$. The condition is what allows $L''$ at the intermediate point of the expansion below to be replaced by $L''$ at $\theta$. For an iid sample $s_n$ the log likelihood is $L(\theta\mid s_n)=\sum L(\theta\mid x_i)$ and, assuming the maximum likelihood estimator $\hat\theta_n$ exists and is consistent, the mean-value expansion of $L'(\hat\theta_n\mid s_n)=0$ about $\theta$ gives
\begin{equation}
\sqrt{n}(\hat\theta_n-\theta)=\frac{n^{-1/2}L'(\theta\mid s_n)}{-n^{-1}L''(\theta^*_n\mid s_n)},
\label{eq:mvt}
\end{equation}
with $\theta^*_n$ between $\theta$ and $\hat\theta_n$. The denominator converges in probability to $I_1(\theta)$; the numerator is asymptotically $N(0,I_1(\theta))$ by the central limit theorem applied to the score; and the conclusion is Bahadur's item 1, attributed to Fisher:
\[
\sqrt{n}(\hat\theta_n-\theta)\ \to\ N\big(0,\,1/I_1(\theta)\big)\quad\text{in distribution under } P_\theta .
\]

Immediately after the statement Bahadur adds a note, and the note is the hinge of this section:
\begin{quote}
This does not assert that $E_\theta(\hat\theta_n)=\theta+o(1)$ or that $\operatorname{Var}_\theta(\hat\theta_n)=\dfrac{1}{nI_1(\theta)}+o(1/n)$.
\end{quote}
The theorem is about the sampling distribution of $\hat\theta_n$ and says nothing about its moments. Convergence in distribution does not imply convergence of means or variances; the additional condition that would deliver the moments---uniform integrability of $n(\hat\theta_n-\theta)^2$---is not part of the theorem and is often false. The binomial supplies the simplest instance: with $\hat p=y/n$, the studentized statistic $(\hat p-p)/\sqrt{\hat p(1-\hat p)/n}$ converges in distribution to $N(0,1)$ and has no finite moments at any $n$, because $\hat p\in\{0,1\}$ has positive probability. Yet the limit theorem is what licenses, in practice, the description of $\hat\theta_n$ by a mean and a variance: ``approximately unbiased with variance approximately $1/(nI_1)$'' is the form in which the theorem is used, and it is the form the theorem does not assert.

The generalized reading of \eqref{eq:mvt} was given in \cite[Sec.~6.3]{NewLook}: the expansion says that $\hat\theta_n$ is, to first order, a relabeled standardized score. Everything approximate in the display is in the relabeling---the denominator, and the passage from the score to its root. The score itself needs no approximation. At every $n$ and every $\theta$,
\[
E\big(\gamma^{(1)}_\theta\big)=0,\qquad \operatorname{Var}\big(\gamma^{(1)}_\theta\big)=nI_1(\theta),\qquad \Lambda\big(\gamma^{(1)}\big)=nI_1(\theta),
\]
the last by Theorem~1 of \cite{NewLook}. What is asymptotic about the score is the \emph{shape} of its distribution, which the central limit theorem describes; what is exact is its information. The variance in the display is the variance of a function on $\Theta$ evaluated at $\theta$---a vertical quantity---and it is exact because it is computed under the distribution the argument names. The variance Bahadur's note declines to assert is the variance of a point estimator---a horizontal quantity, the spread of the roots along the parameter axis---and it is unavailable because the point estimator's distribution is only known in the limit. Fisher's claim that information is a finite-sample quantity is literally true of the generalized estimator and only asymptotically true of the point estimator it induces.

Chapter 8 opens, before the regular case, with Example 7: $X_i$ uniform on $(0,\theta)$, where the likelihood is discontinuous, ``the information function is not defined,'' $\hat\theta_n=\max X_i$ has a non-normal limit, and $\hat\theta_n-\theta=O(1/n)$. This is outside the regular chapters of the lectures and outside the generalized theory as developed here, which requires condition (ii)---that $\tau(s,\cdot)$ be continuously differentiable on $\Theta$---in order to define $\Lambda$ at all; we note it and set it aside, as Bahadur does.

\subsection{Lectures 28--29: Hodges and Fisher's bound}
\label{sec:hodges}

Fisher asserted \cite{Fisher1925}, and it was long believed, that if $\sqrt{n}(t_n-\theta)\to N(0,v(\theta))$ for any sequence of estimates $t_n$, then $v(\theta)\ge 1/I_1(\theta)$. Bahadur records Hodges's 1951 counterexample, first published by Le~Cam \cite{LeCam1953}: $X_i$ iid $N(\theta,1)$, $\Theta=\mathbb{R}$, and
\begin{equation}
t_n=\begin{cases}\bar X_n & |\bar X_n|>n^{-1/4}\\ c\bar X_n & |\bar X_n|\le n^{-1/4}\end{cases}
\label{eq:hodges}
\end{equation}
with $-1<c<1$. Then $\sqrt{n}(t_n-\theta)\to N(0,v(\theta))$ for every $\theta$, with $v(\theta)=1$ for $\theta\ne0$ and $v(0)=c^2<1=1/I_1$. Lecture 29 states the property at issue as a condition on a sequence of estimates, which we display as Bahadur does:
\begin{quote}
\emph{Condition $(\ast\ast\ast)$. $\{t_n-\theta\}$ is $AN(0,v(\theta)/n)$ under $\theta$ (for each $\theta$).}
\end{quote}
Here $\{z_n\}$ is $AN(\mu_n,\sigma_n^2)$, in Bahadur's definition, if $P\{(z_n-\mu_n)/\sigma_n\le z\}\to\Phi(z)$ for every $z$: the distribution function of the standardized sequence converges pointwise to the normal one, which, $\Phi$ being continuous, is convergence in distribution of the standardized sequence.
Hodges's $t_n$ satisfies $(\ast\ast\ast)$ with the $v$ above, and Lecture 29 assembles the repair: the set $\{\theta: v(\theta)<1/I_1(\theta)\}$ has Lebesgue measure zero (item 2, Le~Cam/Bahadur \cite{LeCam1953,Bahadur1964}); if $v$ and $I_1$ are continuous the bound holds everywhere (Corollary); the convolution theorem (item 4, H\'ajek \cite{Hajek1970}) shows that a regular limit distribution is $N(0,1/I_1)$ convolved with independent noise, so its variance, if it exists, is at least $1/I_1$. Bahadur's summary: ``At least in the iid case, Fisher's assertion is essentially correct.''

Two of Bahadur's own remarks along the way should be kept in view. The first is that the bound ``should not be confused with the C--R bound, since [asymptotic normality] does not imply that $t_n$ is unbiased, nor that $v(\theta)\approx n\operatorname{Var}_\theta(t_n)$.'' The variance of the limit is again not the limit of the variances. The second is the form of the repair itself: the classical theory recovers Fisher's bound by excluding the offending estimators (measure zero; regularity) and by proving that among the survivors the normal limit is best (convolution). These are the asymptotic counterparts of the apparatus catalogued in \cite[Sec.~2.4]{NewLook}---unbiasedness, admissibility, minimaxity---and they answer the same need: a pointwise criterion, here the asymptotic variance, cannot by itself single out the maximum likelihood estimator, so the class of competitors must be trimmed until it can.

The generalized theory has nothing to repair, because superefficiency cannot occur. For any generalized estimator $\tau$, at every $\theta$ and every $n$, $\Lambda(\tau)\le nI_1(\theta)$ (\cite[Thm.~1]{NewLook}), with no regularity condition on $\tau$ beyond membership in the class. Hodges's estimator induces two generalized estimators, the score-map image $\tau_{t_n}$ and the expectation-map image $\tau_{E(t_n)}=t_n-E(t_n)$, and the bound applies to both. The score map is the more principled of the two and gives the shorter answer (Section~\ref{sec:hodges-finite}); we take the expectation map first because it is built from the mean and variance of $t_n$, which are the quantities the classical comparison is conducted on, and so it meets the classical account on its own ground. What the asymptotic variance $v(0)=c^2$ records is not that $t_n$ is more informative than $\bar X_n$ at the origin but that its mean function has slope $c$ there: for $\theta$ near $0$ and $n$ large, $E_\theta(t_n)\approx c\theta$, so that
\begin{equation}
\Lambda\big(\tau_{E(t_n)}\big)\Big|_{\theta=0}=\frac{\big(\tfrac{d}{d\theta}E_\theta t_n\big)^2}{\operatorname{Var}_0(t_n)}\approx\frac{c^2}{c^2/n}=n=nI_1 .
\label{eq:hodges-slope}
\end{equation}
The ``super'' in superefficiency is a variance reported without the slope that belongs with it; the information criterion, which divides the one by the other, sees an estimator that is fully efficient at $\theta=0$ and nowhere better than the score.

\subsection{Hodges's estimator at finite $n$}
\label{sec:hodges-finite}

\begin{figure}[t]
\centering
\includegraphics[width=\textwidth]{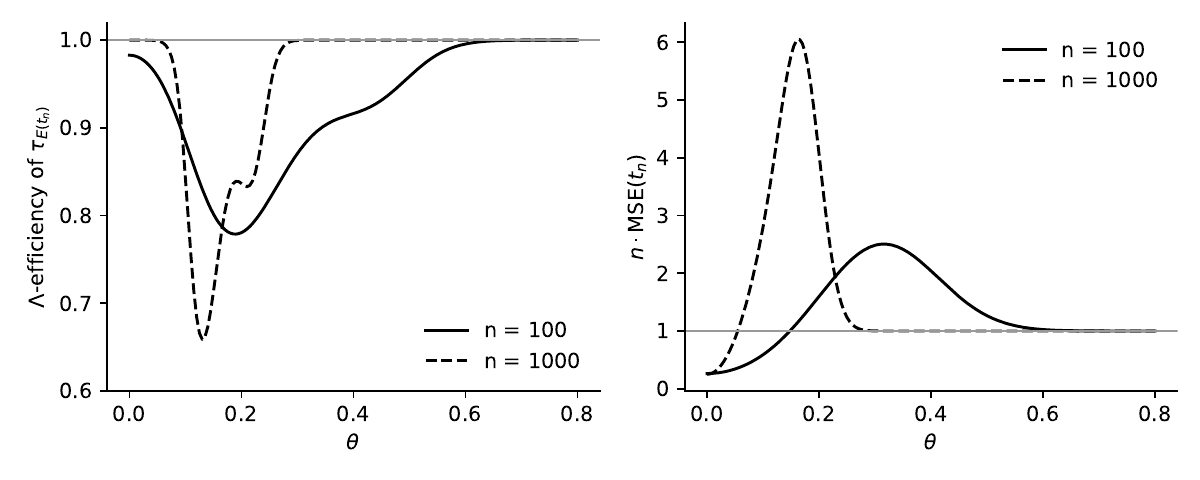}
\caption{Hodges's estimator \eqref{eq:hodges} with $c=1/2$. Left: $\Lambda$-efficiency of the expectation-map image $\tau_{E(t_n)}$ as a function of $\theta$. Right: $n$ times the mean square error. Solid, $n=100$; dashed, $n=1000$. The horizontal line is the value $1$ attained by $\bar X_n$ under both criteria.}
\label{fig:hodges}
\end{figure}

Figure~\ref{fig:hodges} replaces the approximation in \eqref{eq:hodges-slope} by exact computation and extends it over the parameter. Its left panel is the $\Lambda$-efficiency $\Lambda(\tau_{E(t_n)})/(nI_1)$ of Hodges's estimator with $c=1/2$, as a function of $\theta$ for $n=100$ and $n=1000$; its right panel, on the same axis, is $n\cdot\mathrm{MSE}_\theta(t_n)$, the quantity through which Le~Cam's analysis of superefficiency is usually presented. Both panels are exact: the moments of $t_n$ are truncated-normal moments in closed form, recorded at the end of this subsection.

Three features of the figure carry the argument. First, the efficiency never exceeds $1$: at $\theta=0$ it is $0.98$ for $n=100$ and $1.00$ for $n=1000$, in agreement with the slope calculation, and the gray line is a ceiling the curves touch but do not cross. Second, the efficiency dips in the transition zone $|\theta|\approx n^{-1/4}$, where $t_n$ switches between its two branches and its distribution is bimodal: the minimum is $0.78$ near $\theta=0.19$ for $n=100$ and $0.66$ near $\theta=0.13$ for $n=1000$. This is the region in which, on the right, $n\cdot\mathrm{MSE}$ rises to $2.5$ and $6.0$---the phenomenon, due to Le~Cam, that at any fixed $n$ the superefficient estimator is worse than $\bar X_n$ on a neighborhood of the origin whose width shrinks with $n$ but whose depth grows. The information criterion sees the same defect, at the same place, and reports it in the same units as everything else it reports: as a fraction of the Fisher information. Third, outside the transition zone both curves return to $1$, and they do so at every $n$; nothing in the figure requires a limit.

The figure uses the expectation map $\tau_{E(t_n)}$ because that is the map on which the classical comparison---means and variances of $t_n$---is conducted. Under the score map $\tau_{t_n}$, the score of the marginal family of the statistic $t_n$, the conclusion is stronger and needs no figure. For $0<c<1$ the map $\bar x\mapsto t_n(\bar x)$ is injective, since the range $\{c\bar x:|\bar x|\le n^{-1/4}\}=[-cn^{-1/4},cn^{-1/4}]$ and the range $\{\bar x:|\bar x|>n^{-1/4}\}$ are disjoint; $t_n$ therefore generates the same $\sigma$-field as $\bar X_n$, is sufficient, and its score-map image is the full-data score: $I_{t_n}=\Lambda(\tau_{t_n})=nI_1$ at every $\theta$ and every $n$. Under the score map, Hodges's estimator discards no information at all. Any \emph{family-aware} assessment---one that evaluates an estimator at $\theta$ by how it behaves as the distribution moves within the family near $\theta$, rather than by its sampling distribution under $P_\theta$ alone \cite[Sec.~2.4]{NewLook}---must give $t_n$ and $\bar X_n$ the same value, because the two statistics carry the same information about every distribution in the family; what distinguishes them is the mean function, which enters the expectation map through its slope and enters the classical assessment not at all. Superefficiency is what results from assessing a point estimator by the variance of its horizontal distribution alone---a criterion that is not family-aware, and that mistakes a small slope for a small error.

For the record, the exact computation behind Figure~\ref{fig:hodges} is as follows. Since $t_n$ is piecewise linear in $\bar X_n\sim N(\theta,1/n)$, with $a=n^{-1/4}$, $s=n^{-1/2}$, $\alpha=(-a-\theta)/s$, $\beta=(a-\theta)/s$ and $P=\Phi(\beta)-\Phi(\alpha)$,
\[
E_\theta(t_n)=\theta-(1-c)M_1,\qquad E_\theta(t_n^2)=\theta^2+s^2-(1-c^2)M_2,
\]
where $M_1=\theta P-s\{\phi(\beta)-\phi(\alpha)\}$ and $M_2=\theta^2P+2\theta s\{\phi(\alpha)-\phi(\beta)\}+s^2\{P+\alpha\phi(\alpha)-\beta\phi(\beta)\}$ are the moments of $\bar X_n$ restricted to $|\bar X_n|\le a$, and $\tfrac{d}{d\theta}E_\theta(t_n)=1-(1-c)\{P-\theta(\phi(\beta)-\phi(\alpha))/s-\beta\phi(\beta)+\alpha\phi(\alpha)\}$.

\subsection{Lectures 29--30: the bound is proved by a test, and the LAN condition is a statement about the score}
\label{sec:L29}

Bahadur observes that condition $(\ast\ast\ast)$ implies $t_n$ is asymptotically median unbiased, $P_\theta(t_n\le\theta)\to\tfrac12$, and it is this consequence, assumed to hold uniformly in a neighborhood, that his proof of the bound (item 3) actually uses. The proof deserves attention for its method. Fix $\theta$ and $a>0$, put $\theta_n=\theta+a/\sqrt{n}$, and consider testing $\theta$ against $\theta_n$. The optimal test statistic is the likelihood ratio $\Omega_{\theta_n,\theta}(s_n)$, whose logarithm is, up to negligible terms,
\[
\log\Omega_{\theta_n,\theta}(s_n)=L(\theta_n\mid s_n)-L(\theta\mid s_n)=a\sqrt{I_1(\theta)}\,\bar\gamma^{(1)}_\theta(s_n)-\tfrac12 a^2I_1(\theta),\qquad \bar\gamma^{(1)}_\theta(s_n)=\frac{\gamma^{(1)}_\theta(s_n)}{\sqrt{nI_1(\theta)}} ,
\]
so that the standardized score $\bar\gamma^{(1)}_\theta(s_n)$, which Bahadur writes $K_n$, is equivalent to the likelihood-ratio statistic and is asymptotically $N(0,1)$ under $P_\theta$ and $N(a\sqrt{I_1},1)$ under $P_{\theta_n}$. The Neyman--Pearson lemma \cite{NeymanPearson1933,LehmannRomano2005} then compares the critical region $\{t_n\ge\theta_n\}$, which has asymptotic size at least $\tfrac12$ under $\theta_n$ by median unbiasedness, with $\{\bar\gamma^{(1)}_\theta(s_n)\ge z\}$, and yields $v(\theta)\ge1/I_1(\theta)$.

Read with the definitions recalled above, the argument is Fisher's continuum of significance tests, evaluated at the pair of hypotheses $\theta$ and $\theta+a/\sqrt{n}$, with the standardized score as the test statistic. The classical theory arrives at the bound on the asymptotic variance of a point estimator by way of the generalized estimator and then hands the conclusion back to the point estimator. The information $\Lambda$ is the same comparison made at fixed $n$ and in the limit $a\to0$: the numerator $E(\tau')$ is the rate at which the test statistic responds as the alternative moves away from the null, and the denominator standardizes it. Bahadur's proof needs the alternative at distance $a/\sqrt{n}$ because the point estimator must be given a distribution before it can be compared with anything, and only in the limit does it have one; the score has one at every $n$.

Lecture 30 replaces the iid assumption by the local asymptotic normality condition of Le~Cam \cite{LeCam1960},
\[
\log\Omega_{\theta+a/\sqrt{n},\theta}(s_n)=a'z_n(\theta)-\tfrac12a'I_1(\theta)a+\Delta_n(\theta,s_n),
\]
with $z_n\to N(0,I_1)$ and $\Delta_n\to0$ in probability, under which Le~Cam's and H\'ajek's theorems and the Newton-step results of Chapter 7 (Lectures 25--27, restated in Lecture 30 as item $1^p$) go through unchanged. In the iid case $z_n(\theta)$ is $n^{-1/2}\gamma^{(1)}_\theta(s_n)$: the condition asserts that the log likelihood ratio at local alternatives is, to first order, linear in the score with a quadratic correction of size $I_1$, and that the score is asymptotically normal. LAN is a condition on the score, not on any estimator. Bahadur's Note ii on the condition is the observation that matters here: ``the right-hand side in LAN with $\Delta_n$ omitted is exactly the log-likelihood in the multivariate normal translation-parameter case.'' That is, the asymptotic theory works by showing that every regular family looks, locally and in the limit, like a normal location family. The normal location family is an exponential family, the case in which the score is a linear function of a sufficient statistic and attainment of the information bound is exact \cite[Sec.~5.3]{NewLook}; it is moreover the self-dual exponential family, whose geometry is Euclidean rather than the dually flat geometry that a general exponential family carries \cite{Efron1975,Amari1985,AmariNagaoka2000}. The limit replaces the family at hand by the simplest one there is. The generalized theory does not need the family to look flat. It needs $\Lambda$, which is defined for every generalized estimator in every regular family at every $n$, and which reduces to the classical quantities exactly where the classical theory can compute them.

Chapter 8 closes with the sentence: ``for large $n$, the $N(0,I_1^{-1}(\theta)/n)$ distribution is nearly the best possible for estimates of $\theta$.'' Every qualification in it---``large $n$,'' ``nearly,'' ``estimates of $\theta$''---is the cost of asking a point in $\Theta$ to carry what a function on $\Theta$ carries. For the score the corresponding sentence has no qualifications: at every $n$, $\Lambda(\gamma^{(1)})=nI_1$, and no generalized estimator does better.

\section{Fisher information for a continuum of tests}
\label{sec:fisher}

The previous section used the quantity $\Lambda$ as \cite{VosWu,NewLook} defined it. This section argues that $\Lambda$ is Fisher information. Three quantities carry Fisher's name here, and it helps to state them together. We write a prime for differentiation in the parameter, as Bahadur does.

The \emph{Fisher metric} is the variance of the score $\gamma^{(1)}$: $I_M=\operatorname{Var}(\gamma^{(1)})$. It describes the local behaviour of the distributions in the family $M$.

The \emph{Fisher information in a statistic} $T$ is the negative of the mean derivative of the score $\gamma^{(1)}_T$ of $T$: $I_T=-E\big(\gamma^{(1)\prime}_T\big)$. It describes the information in $T$ for distinguishing between distributions in $M$; Section~\ref{sec:two} argues that its proper form is the standardized one, $I_T=\big(E\,\bar\gamma^{(1)\prime}_T\big)^2$, which has the same value; here $\bar\gamma^{(1)}_T=\gamma^{(1)}_T/\sqrt{\operatorname{Var}(\gamma^{(1)}_T)}$ is the standardized score of $T$.

The \emph{Fisher information in a continuum of tests} $\tau$ is the square of the mean derivative of the standardized test: $I_\tau=\big(E\,\bar\tau'\big)^2$, where $\bar\tau=\tau/\sqrt{\operatorname{Var}(\tau)}$.

Each definition extends the one before it: from the family, to inference based on a statistic, to inference based on a continuum of tests. When $\gamma^{(1)}_T=\gamma^{(1)}$ the values of the Fisher metric and of the information in $T$ are equal, $I_T=I_M$; when $\tau=\gamma^{(1)}_T$ the values of the Fisher information in $T$ and in $\tau$ are equal, $I_\tau=I_T$. The extension from $I_M$ to $I_\tau$ rests on statistical reasoning, set out in Section~\ref{sec:two}; the motivation for $\Lambda$ in \cite{VosWu} was geometric, and is recalled in Section~\ref{sec:metric}. Since $\Lambda(\tau)=I_\tau$, $\Lambda(\tau)$ should be called the Fisher information in $\tau$. The reason for insisting on the name is that ``information'' is used in a general sense and in several technical ones; the quantity defined here is information in the sense of Fisher and in no other, and $\Lambda$-information should not be thought of as something different from Fisher information; it is an extension of Fisher information in a statistic to Fisher's notion of a continuum of significance tests. Sections~\ref{sec:t}--\ref{sec:asy} then draw the consequences for a test in the presence of a nuisance parameter, for variance, and for asymptotics.

\subsection{Fisher information in a statistic: derivative and standardization}
\label{sec:two}

Let $T$ be a statistic whose marginal family $\{P_\theta\circ T^{-1}\}$ is regular, with score $\gamma^{(1)}_T$. The conventional definition takes the Fisher information in $T$ to be the variance of the score, $E\big[(\gamma^{(1)}_T)^2\big]$, and presents $-E\big(\gamma^{(1)\prime}_T\big)$ as a computational convenience, the two being equal by the Bartlett identity \cite{Bartlett1953} (Bahadur's item 13). The equality is of values, not of meanings, and the conventional usage is correspondingly ambiguous: the observed information is defined through the derivative, as $-\gamma^{(1)\prime}_T$, a random variable whose expectation is $I_T$ \cite{EfronHinkley1978}, while the expected information is defined through the variance. The variance is a property of the vertical distribution of the score at a single $\theta$; the mean derivative is a comparison of $\theta$ with its neighbours. The derivative is the essential ingredient, because it is the derivative that makes the quantity an information: a statistic whose score did not respond to changes in the null hypothesis would carry no information about it however small its variance. Defining Fisher information as the variance reinforces the reading of information as inverse variance and leaves unexplained why a variance should be called an information at all.

The derivative alone is not enough, because a rate of response can be made as large as one likes by rescaling the statistic and must be put in standard-deviation units before it means anything. Standardization is the ordinary statistical principle for comparing quantities measured on different scales, and here it is also what the geometry requires, as Section~\ref{sec:metric} shows. For the score the standardized and unstandardized derivatives are related by a simple identity. With $\bar\gamma^{(1)}_T=\gamma^{(1)}_T/\sqrt{I_T}$, and because the expectation of the score is zero,
\[
E\big(\bar\gamma^{(1)\prime}_T\big)=\frac{E\big(\gamma^{(1)\prime}_T\big)}{\sqrt{I_T}}=-\sqrt{I_T},\qquad\text{so that}\qquad I_T=\big(E\,\bar\gamma^{(1)\prime}_T\big)^2 .
\]
The two forms have the same value. They do not have the same meaning, and units of measurement show that the standardized form gives the more interpretable definition.

Location and scale parameters usually carry units; for example, a mean diastolic blood pressure is in mmHg. Statistical significance is parameter-invariant, but clinical or scientific significance is ordinarily stated in such units, so they cannot be ignored when a quantity is to be interpreted. The likelihood and the log likelihood, through which inference is conducted, are unitless as functions on the family: a likelihood ratio is a pure number. The score is not. If $\theta$ is measured in a unit $u$, the score, a derivative of a unitless quantity with respect to $\theta$, has units $u^{-1}$; its derivative has units $u^{-2}$, and so does its variance. The unstandardized form $-E\big(\gamma^{(1)\prime}_T\big)$ is therefore the mean rate of change, per unit $u$, of a quantity that itself has units $u^{-1}$, and its square root, which has the units one wants, is not the rate of change of anything. The standardized score is unitless, as the log likelihood is, and $E\big(\bar\gamma^{(1)\prime}_T\big)$ is the mean rate of change of a unitless test statistic per unit of $\theta$: the number of standard deviations by which the statistic moves when the hypothesized blood pressure moves by one mmHg. This is $\sqrt{I_T}$, obtained directly as a rate, and the information is its square. The situation is the familiar one of the standard deviation and the variance. The spread of a population, or of the sampling distribution of an estimator, is described in the units of the measurement by the standard deviation; the variance, in squared units, is the quantity that adds over independent components. In the same way $\sqrt{I_T}$, a slope in standard deviations per unit, is the quantity one interprets, and $I_T$ is the quantity that adds over independent observations, for the geometric reason given in Section~\ref{sec:metric}. To use the score for inference as the likelihood is used, it should be made unitless, and the definition of Fisher information that does this is the standardized one. We therefore prefer $I_T=\big(E\,\bar\gamma^{(1)\prime}_T\big)^2$ as the definition, with the unstandardized mean derivative and the variance of the score as two expressions that agree with it in value.

Standardization also removes a dependence on the parameterization that the score carries and that has nothing to do with inference. Inference uses points of the sample space to distinguish among the distributions of $M$; a parameter $\theta$ labels those distributions but need not capture the properties that matter for distinguishing them. The value of the score at a given distribution depends on the label---under a reparameterization $\xi$, $\gamma^{(1)}_\xi=\gamma^{(1)}_\theta\,d\theta/d\xi$---while the value of the standardized score does not. Binomial and negative binomial sampling make the point concrete \cite{VosSLP}. Observing $y$ successes in $n$ Bernoulli trials, with $n$ fixed in the first plan and $y$ fixed in the second, gives likelihood functions proportional to each other and therefore identical scores, $(y-np)/\{p(1-p)\}$, as functions of $p$. The two families are nonetheless different, with Fisher metrics $n/\{p(1-p)\}$ and $y/\{p^2(1-p)\}$, and the standardized scores are different:
\[
\bar\gamma_{\mathrm{Bin}}(p)=\frac{y-np}{\sqrt{np(1-p)}},\qquad \bar\gamma_{\mathrm{NB}}(p)=\frac{y-np}{\sqrt{y(1-p)}} .
\]
For $y=6$, $n=20$ and the null value $p=0.5$, the first is $-1.79$ and the second $-2.31$. The score shows only the relationship between the sample point and values in the parameter space, and that relationship is the same for the two plans; what is relevant for inference is how the test responds as the distribution moves within the family, and the standardized score, which is computed in the family the sampling plan actually generates, reflects the difference between the plans. The bearing of this on the strong and weak likelihood principles is the subject of \cite{VosSLP}.

Once the Fisher information in a statistic is written through the standardized score, the extension from the score to any smooth test is nearly immediate: replace the score by a generalized estimator. For a generalized estimator $\tau$ with standardization $\bar\tau$,
\begin{equation}
I_\tau=\big(E\,\bar\tau'\big)^2=\frac{\big(E\,\tau'\big)^2}{E(\tau^2)},
\label{eq:lambda}
\end{equation}
the second equality again because $E\tau=0$. This is $\Lambda(\tau)$. For $\tau=\gamma^{(1)}$ the Bartlett identity gives $-E\,\gamma^{(1)\prime}=E\big[(\gamma^{(1)})^2\big]$, and \eqref{eq:lambda} collapses to either one. The collapse is a convenience for the score and a trap for the concept: it permits Fisher information to be defined by the variance alone, and once it is so defined the question of what Fisher information means for a continuum of tests other than the score cannot be asked, because the definition has discarded the ingredient that would answer it. Fisher, having shown that the score is the most informative test in the continuum, had no occasion to write the general form.

One consequence of the extension should be noted. A statistic $T$ supports more than one continuum of tests. Its score continuum, the score-map image $\tau_T$ of Section~\ref{sec:ch8}, has $I_{\tau_T}=I_T$ at every $n$; its expectation-map image $\tau_{E(T)}=T-E(T)$ can have less, the two agreeing exactly when $T$ is the natural statistic of its own marginal exponential family \cite[Prop.~1]{MSE}. Fisher information therefore attaches to a continuum of tests, and the Fisher information in a statistic is the Fisher information in its score continuum. We write $I_\tau$ when a particular continuum is assessed and keep $\Lambda(\cdot)$ for statements, such as the bound below, that concern Fisher information as a function on the space of generalized estimators.

\subsection{The Fisher metric, and the Fisher upper bound}
\label{sec:metric}

That the variance of the score describes the family itself was noticed by Rao \cite{Rao1945}: it is a Riemannian metric on $M$, describing the local structure of the family at each distribution, and locally it is twice the Kullback--Leibler divergence \cite[Sec.~2.2]{MSE}. $I_M$ is the $I$ of Section~\ref{sec:ch8}, the subscript recording that it belongs to the family rather than to any parameterization of it. The Fisher metric, being a smooth inner product on the tangent spaces of $M$, is conceptually distinct from the Fisher information in a statistic. The derivative that the information in a statistic requires is not an object of the tangent space at all: the derivative of a tangent vector need not be a tangent vector, and in a curved family the derivative of the score leaves the tangent space, the departure being Efron's statistical curvature \cite{Efron1975}.

Where the derivative goes is the Hilbert space. A statistical manifold---a family of distributions with common support, carrying the Fisher metric---is a richer structure than a Riemannian manifold alone: at each $m\in M$ it carries the Hilbert space $V_m$ of square-integrable functions on the sample space, which contains the tangent space $T_mM$ as the span of the score \cite{VosWu}. The derivative of a tangent vector can leave $T_mM$ but stays in $V_m$, and so do the generalized estimators. Two geometries are attached to every continuum of tests, the score first among them. For a fixed sample point $s$, the curve $\tau_s$ orders the distributions of $M$ by their inconsistency with $s$; for a scalar parameter it can be drawn in $\Theta\times\mathbb{R}$, but what is relevant is how it changes as the distribution moves in $M$, not the parameter values along the axis. For a fixed distribution $m_\circ$, the vertical slice $\tau(\cdot,m_\circ)$ orders the points of the sample space by their extremity under $m_\circ$; since a test statistic and any positive multiple of it give the same ordering, the slice is a direction in $V_{m_\circ}$, the unit vector $\bar\tau(\cdot,m_\circ)$ orthogonal to the constants. The two orderings point the same way: a sample point that is extreme under a model is a sample point with which the model is inconsistent. For the score the first geometry gives the score curve, and the second gives a direction that lies in the tangent space. Replacing the score by another smooth test replaces the score curve by another smooth curve and the score direction by another direction in $V_m$. The structure of a statistical manifold thus provides, in $V_m$, a large class of tests containing the score as the special case whose direction lies in $T_mM$.

The main result of \cite{VosWu}, Theorem~1 of \cite{NewLook}, can now be written with every ingredient visible. We propose to call it the \emph{Fisher upper bound}:
\begin{equation}
\frac{\big(E\,\tau'\big)^2}{E(\tau^2)}\ \le\ I_M ,
\label{eq:fiub}
\end{equation}
with equality if and only if $\tau(\cdot,\theta)$ is almost surely a multiple of the score. On the left is the rate at which the continuum of tests $\tau$ responds to movement in $M$, standardized by its variance; on the right is an intrinsic property of $M$. The prime is written, rather than a derivative in a named parameter, so that no parameter appears on one side only: a change of parameter from $\theta$ to $\xi$ multiplies both sides by the same constant, $(d\theta/d\xi)^2$, so the inequality, and the ratio of its two sides, are the same in every parameterization. The exact form $I_\tau=\mathrm{Corr}^2(\tau,\gamma^{(1)})\,I_M$ gives the bound its geometric reading. $\mathrm{Corr}(\tau,\gamma^{(1)})$ is the cosine of the angle between the direction $\bar\tau$ and the tangent space, so the optimal test is one whose direction lies in the tangent space, and the suboptimality of any other test is measured by the departure of its direction from it. It is not a coincidence that Fisher's construction and Rao's metric agree at the optimum: a continuum of tests can respond to movement in the family no faster than the family itself moves, and the rate at which the family moves is the metric. The same reading explains why higher-order derivatives of the likelihood, which enter the classical theory through the Bhattacharyya bounds \cite{Bhattacharyya1946}, cannot improve a test \cite[Sec.~4.2]{NewLook}: their components orthogonal to the tangent space add variance to a test and never information.

The Fisher upper bound should be set beside the Cram\'er--Rao lower bound \cite{Rao1945,Cramer1946}, $\operatorname{Var}(t)\ge (g')^2/I$ for $t$ unbiased for $g(\theta)$, from which it differs algebraically only by a rearrangement \cite[eq.~(12)]{NewLook}. The geometry shows the difference to be fundamental. The Cram\'er--Rao bound places a quantity that is not family-aware---the variance of $t$, an attribute of the horizontal distribution of $t$ under the single distribution $P_\theta$---against a bound that depends on the family; the two sides are objects of different kinds, and the estimator being bounded changes with the parameter, since an unbiased estimator of $g(\theta)$ and an unbiased estimator of a non-affine transformation of $g(\theta)$ are different statistics, each with its own bound. The Fisher upper bound places a family-aware assessment against a family bound. Both sides of \eqref{eq:fiub} are built from the score: the right side is its squared length and the left side, since $E\,\bar\tau'=-\langle\bar\tau,\gamma^{(1)}\rangle$, is the squared length of its projection on the direction $\bar\tau$. Both are $(0,2)$-tensors on $M$, which is why a change of parameter multiplies them by the same constant; the same $\tau$ is assessed against the same bound whatever labels are used, and the efficiency $\mathrm{Eff}(\tau)=I_\tau/I_M$ is a scalar on $M$. This is what the notation $I_M$, rather than $I_\Theta$, is meant to record. The geometry also explains why the mean slope enters squared. The proof of \eqref{eq:fiub} is a projection: the score decomposes as its component along $\bar\tau$ plus an orthogonal remainder, and the Pythagorean relation in $V_\theta$ makes squared lengths additive over orthogonal components. Fisher information is the squared length of a projection because that is the scale on which orthogonal contributions add, and it is on that scale that information from independent observations is additive and that the efficiency of a statistic is the fraction of information it retains.

Given the optimality of the score, one may ask why other tests need be considered at all. There are two reasons. Criteria other than local optimality can matter, robustness to departures from the model being the most familiar; a test chosen for robustness is a direction in $V_m$ away from the tangent space, and \eqref{eq:fiub} prices the departure in Fisher's units. And the score may not be available: when a nuisance parameter is present the score for the parameter of interest depends on the nuisance parameter, and the test one can actually perform is a different continuum, as Section~\ref{sec:t} shows for Student's $t$.

\subsection{A test can lose Fisher information that its statistic contains}
\label{sec:t}

Let $X_1,\dots,X_n$ be iid $N(\mu,\sigma^2)$ with both parameters unknown and $\mu$ the parameter of interest, and let $\bar X$ and $S$ be the sample mean and the sample standard deviation, $S^2=\sum(X_i-\bar X)^2/(n-1)$. (In this subsection $S$ is the sample standard deviation, as is customary, and not the sample space.) The statistic is the pair $(\bar X,S)$. The score for $\mu$ is $n(\bar X-\mu)/\sigma^2$, which is already orthogonal to the score for $\sigma$ (the two parameters are orthogonal in the normal family), so the nuisance-orthogonalized Fisher information for $\mu$ in $(\bar X,S)$ \cite[Sec.~3]{VosWu} is $I^\perp=n/\sigma^2$. All of it is carried by $\bar X$: the distribution of $S$ is free of $\mu$ and $S$ is independent of $\bar X$, so $S$ provides no information for $\mu$. The score continuum, standardized, is the $z$-test, $\sqrt n(\bar X-\mu)/\sigma$; it has Fisher information $I^\perp$, but it depends on $\sigma$ and cannot be computed when $\sigma$ is unknown.

The test one actually performs when $\sigma$ is unknown is Student's $t$-test \cite{Student1908}, which uses $S$ in place of $\sigma$. As a continuum of tests it is a function of the statistic $(\bar X,S)$ and of $\mu$ alone,
\[
\tau(\mu)=\frac{\bar X-\mu}{S} .
\]
This is a generalized estimator for $\mu$ in the sense of \cite[Def.~5]{VosWu}: it has mean zero at every $(\mu,\sigma)$, it is smooth in $\mu$, its variance $E[(\bar X-\mu)^2]\,E[S^{-2}]=(\sigma^2/n)E[S^{-2}]$ is finite for $n\ge4$, and it is orthogonal to the nuisance score $\partial\ell/\partial\sigma=-n/\sigma+\sum(X_i-\mu)^2/\sigma^3$ because $\bar X-\mu$ is independent of $S$ and symmetric about zero, so that every term in $E[\tau\,\partial\ell/\partial\sigma]$ vanishes. Its Fisher information is computed from \eqref{eq:lambda} directly: $E[\partial\tau/\partial\mu]=-E[S^{-1}]$, and
\[
I_\tau=\frac{\big(E[S^{-1}]\big)^2}{(\sigma^2/n)\,E[S^{-2}]}=\frac{n}{\sigma^2}\cdot\frac{\big(E[S^{-1}]\big)^2}{\sigma^2E[S^{-2}]} .
\]
With $\nu=n-1$, $\nu S^2/\sigma^2\sim\chi^2_\nu$, so $E[S^{-2}]=\nu/\{(\nu-2)\sigma^2\}$ and $E[S^{-1}]=\sqrt{\nu}\,\Gamma\big(\tfrac{\nu-1}{2}\big)/\{\sqrt2\,\Gamma(\tfrac{\nu}{2})\,\sigma\}$, and the efficiency of the $t$-test continuum relative to $I^\perp$ is
\begin{equation}
\mathrm{Eff}(\tau)=\frac{I_\tau}{I^\perp}=\frac{(\nu-2)\,\Gamma\!\big(\tfrac{\nu-1}{2}\big)^2}{2\,\Gamma\!\big(\tfrac{\nu}{2}\big)^2},\qquad \nu=n-1\ge3 ,
\label{eq:teff}
\end{equation}
free of $\mu$ and $\sigma$. Its values are $2/\pi=0.637$ at $n=4$, $\pi/4=0.785$ at $n=5$, $0.931$ at $n=10$, $0.971$ at $n=20$, $0.989$ at $n=50$, and $0.995$ at $n=100$; the deficit $1-\mathrm{Eff}(\tau)$ is close to $1/(2n)$ for large $n$. The $t$-test, as a continuum of tests, uses less Fisher information for $\mu$ than its statistic $(\bar X,S)$ contains, by an amount equivalent to about half an observation at every sample size. $S$ carries no information for $\mu$, but it is what frees the test from $\sigma$, and the loss is the price.

The two geometries of Section~\ref{sec:metric} show where the loss comes from. For the $z$-test the standardized curves over the $\mu$-axis are parallel straight lines, one for each sample, all with slope $-\sqrt n/\sigma$. For the $t$-test the curves are again straight lines, but the slope of the standardized curve, $\bar\tau'=-S^{-1}/\sqrt{\operatorname{Var}(\tau)}$, varies from sample to sample. The squared length of $\bar\tau'$ in the Hilbert space is $E\big[(\bar\tau')^2\big]=E[S^{-2}]/\{(\sigma^2/n)E[S^{-2}]\}=n/\sigma^2$, the whole of $I^\perp$, and it divides into the squared length of the projection of $\bar\tau'$ on the constant functions and the variance that remains:
\[
I^\perp=\big(E\,\bar\tau'\big)^2+\operatorname{Var}(\bar\tau')=I_\tau+\operatorname{Var}(\bar\tau') .
\]
The Fisher information lost by the $t$-test is exactly the variance of the slopes of its curves, and $\mathrm{Eff}(\tau)=(E[S^{-1}])^2/E[S^{-2}]$. This is the remark of Section~\ref{sec:metric} about variance without information, seen in the curves: the $t$ statistic is the score multiplied by a random factor independent of it, and the variability of that factor is the part of the test's direction that lies outside the tangent space.

Two remarks. First, the loss is real and it is finite-sample: it is the price of standardizing by an estimated scale, and first-order asymptotics cannot see it, since $\mathrm{Eff}(\tau)\to1$. It is the same phenomenon as the binomial test with estimated standard error of Section~\ref{sec:L28}, whose studentized statistic is undefined on an event of positive probability at every $n$, so that it lies outside the Hilbert space and the loss is total; here $S$ is bounded away from zero often enough for the loss to be finite. The boundary between the two is visible in \eqref{eq:teff}: for $n\le3$, $E[S^{-2}]=\infty$ and the $t$-test continuum is not in the Hilbert space at all. Second, the loss is local, and the phrase ``half an observation'' should be read accordingly. $I_\tau$ is the squared mean slope of the standardized curves $\bar\tau_s$ over $\Theta$, computed at each $\theta$ from the derivative there; it describes how the vertical distributions of the $t$-test continuum turn as $\theta$ moves, and nothing more. Fisher information is additive in observations, so a deficit of $1-\mathrm{Eff}(\tau)$ can be quoted as a number of observations, and that is the sense of the phrase. It does not say that the $t$-test with $n$ observations is equivalent to the $z$-test with $n-\tfrac12$: the two continua have vertical distributions of different shapes ($t_{n-1}$ and normal), and their tail areas differ in more than slope. The same caution was drawn in \cite[Sec.~4.1]{VosWu}, where efficiency translated into sample size approximated one estimator by another well for normal data and poorly for $t_3$ data.

\subsection{Why Fisher information, not variance, carries the limit theorem}
\label{sec:why}

Section~\ref{sec:L28} observed that the central limit theorem, as Bahadur uses it, describes the distribution of the score and says nothing about the moments of the estimator built from it. The two ingredients of \eqref{eq:lambda} explain why that theorem and Fisher information sit together comfortably while the theorem and variance do not. The theorem is a statement about $\gamma^{(1)}_\theta$, the score evaluated at the same $\theta$ that supplies the distribution: not a true distribution, but a hypothesized one, under which one asks how extreme the sample would be. Both the mean derivative and the variance in \eqref{eq:lambda} are computed on that same object under that same $\theta$, and both are exact. The variance of $\hat\theta_n$ is computed on a different object, whose distribution the theorem approximates and whose moments it does not control. Fisher information is defined on the object the theorem is about; variance is defined on the object the theorem is used for. This is the sense in which variance is the wrong quantity to describe the limit theorem and Fisher information the right one---not because the numerical values differ in the limit, for they need not, but because one of them is exact at every $n$ and the other is not.

\subsection{What asymptotics is for}
\label{sec:asy}

An asymptotic result enters inference as a justification for a procedure applied to the $n$ at hand, and a justification is only as strong as the link between the limit and that $n$.

For an assessment by variance the link is weak in two places. The limit of the finite-sample moments need not be the moments of the limiting distribution: this is Bahadur's note after item~1, and the studentized binomial statistic of Section~\ref{sec:L28}, which has a normal limit and no finite moments at any $n$, is the simplest instance. And where the moments do converge, the variance at the $n$ in hand is still not the variance the limit reports, by an amount the limit does not state. The classical theory compensates with further apparatus---uniform integrability where it can be had \cite[Sec.~2.5]{vanderVaart1998}, the regularity conditions and measure-zero exceptions of Lectures 28--29 where it cannot \cite{LeCam1953,Bahadur1964,Hajek1970}---and each piece of apparatus is a further condition that the $n$ in hand must be assumed to satisfy.

For an assessment by Fisher information neither weakness arises, because the assessment is not obtained from a limit. The Fisher upper bound \eqref{eq:fiub} holds at every $n$ for every generalized estimator; the score attains it at every $n$; the efficiency of the $t$-test continuum \eqref{eq:teff} and of Hodges's estimator (Figure~\ref{fig:hodges}) are exact finite-sample quantities. What asymptotics is still needed for is narrower: the shape of a vertical distribution, when the tail area of an observed value must be computed and the exact distribution is unavailable. That is a real use, and the central limit theorem serves it well, because shape is what convergence in distribution delivers. It is a use of asymptotics to approximate, not to assess. The justification of a procedure---that it uses the Fisher information available, or how much of it it loses---does not wait on the limit, and so it is as strong at $n=10$ as at $n=10{,}000$. Asymptotics serves inference at a finite sample size by supplying a shape; Fisher information, not variance, supplies the assessment.

\section{Discussion}
\label{sec:discussion}

The title names three things, and the paper has tried to say what each is for. Fisher information is the assessment of a continuum of significance tests: the squared mean rate at which the standardized continuum responds to the parameter, exact at every sample size, bounded above by an intrinsic property of the family, and attained by the score. Variance has a single role in that assessment, the standardization of the rate; its classical role as the assessment itself belongs to a different object, the point estimator, and to the horizontal distributions that object has and the continuum does not. Asymptotics supplies the shape of a sampling distribution when the exact shape is unavailable (Section~\ref{sec:asy}); it no longer has to assess, and Fisher's claim that information is a small-sample quantity, which the classical theory could defend only asymptotically, is true of the object he was describing.

Reading Bahadur's Chapter 8 in these terms does not overturn any of it. The limit theorems remain theorems, and the estimators they describe remain the estimators one computes. What changes is what the chapter is seen to be about. Each of its three parts reaches an exact statement about the score and then hands the conclusion to a point estimator at the cost of an approximation: the mean-value expansion of Lecture 28, the local-alternative test in the proof of Fisher's bound, the linearity in the score that the LAN condition asserts. Read as statements about the score they are exact and already proved. Read as statements about the estimator they are approximations, and Bahadur's Note after item 1 says how far the approximation reaches. Hodges's estimator is the case in which the two readings come apart most visibly, and the exact computation of Section~\ref{sec:hodges-finite} shows the information criterion seeing, at every $n$ and in the same units as everything else it reports, the defect that the asymptotic account can describe only through Le~Cam's theorem on the shrinking neighborhood.

One of the paper's proposals is terminological, and we think it matters. What \cite{VosWu} called $\Lambda$-information is Fisher information, defined as Fisher's continuum of tests requires and extended from statistics to the continua themselves. The proposal is not a renaming for its own sake. Defining Fisher information as the variance of the score, which the Bartlett identity permits, makes it impossible to ask what information a test other than the score utilizes; defining it by the two ingredients makes the question immediate, and the $t$-test continuum shows that the answer can be exact, finite-sample, and invisible to first-order asymptotics. The bound $\Lambda(\tau)\le I_M$, written out, is a statement that a continuum of tests cannot respond to movement in the family faster than the family moves, and its geometry---vertical distributions, tensors on $M$, the Pythagorean additivity of squared lengths---is what separates it from the Cram\'er--Rao inequality with which it shares an algebraic form.

One consequence of the reading deserves mention, because Bahadur's Chapter 8 supplies an instance of it. Fisher held that a significance test needs no alternative hypothesis \cite{Fisher1955,Fisher1956}, and the position has generally been received as a mistake \cite{NeymanPearson1933,Neyman1956,Lehmann1993} or as one of his obscurities \cite{Savage1976}. If his information is the quantity \eqref{eq:lambda}, the position is a consequence of the definition. The numerator is a derivative in $\theta$ at the null; it is the response of the test to all the local alternatives at once, in the direction the family moves. A test assessed by its information has already been compared with every alternative in a neighborhood of the null. Bahadur's proof of Fisher's bound is the classical theory doing the same thing in slow motion: it names an alternative $\theta+a/\sqrt n$, tests against it with the score, and lets $n\to\infty$, which is to say it lets the alternative approach the null. The information criterion makes the comparison at fixed $n$ and in the limit $a\to0$, and needs no alternative to be named because the derivative has already taken all of them into account. What Fisher declined to supply was a fixed alternative; the local ones are built into the object he was using.

Three limits of the treatment should be stated. The parameter is scalar; the multiparameter theory, with nuisance parameters handled by orthogonalization, is in \cite{VosWu}, and the $t$ example used only the simplest case of it. The families are smooth, as in Bahadur's regular case; Example 7 of his Chapter 8, the uniform on $(0,\theta)$, is outside the present theory as it is outside his, and what plays the role of $\Lambda$ when the score does not exist remains open. And the asymptotic theory of tests, Bahadur's own next subject, has not been touched: his exact slope \cite{Bahadur1960,Bahadur1967,Bahadur1971} assesses a test by a tail area rather than a moment, and its relation to the local, exact assessment by Fisher information is a natural next question.

Bahadur ends Chapter 8 with the sentence ``for large $n$, the $N(0,I_1^{-1}(\theta)/n)$ distribution is nearly the best possible for estimates of $\theta$.'' The reading offered here removes the qualifications one at a time. ``Estimates of $\theta$'' become continua of tests over $\Theta$; ``nearly'' becomes equality, attained by the score; ``for large $n$'' becomes every $n$. What remains is the Fisher upper bound, which is what the chapter was proving all along.

\section*{Acknowledgments}
During the preparation of this manuscript, the author used Claude (Anthropic)
to assist in drafting and revising the text, organizing the exposition,
checking and formatting the references, and preparing the \LaTeX{} source.
The author has reviewed and edited the output, has verified the mathematical
and statistical content, and takes full responsibility for the content of
this publication.

\end{document}